\documentclass[11pt]{amsart}

\usepackage{amsmath,amsxtra,latexsym,amsthm,amssymb,amscd,amsfonts}
\usepackage{graphicx}

\usepackage{color,enumerate}
\usepackage{lscape}

\newtheorem{theorem}{Theorem}[section]
\newtheorem{corollary}[theorem]{Corollary}
\newtheorem{lemma}[theorem]{Lemma}

\theoremstyle{definition}
\newtheorem{definition}[theorem]{Definition}
\newtheorem{example}[theorem]{Example}
\newtheorem{remark}[theorem]{Remark}

\numberwithin{equation}{section}

\def\Limsup{\mathop{{\rm Lim}\,{\rm sup}}}
\def\Liminf{\mathop{{\rm Lim}\,{\rm inf}}}

\def\oz{\bar{z}}

\def\ox{\bar{x}}

\def\gph{\mbox{\rm gph}\,}
\def\epi{\mbox{\rm epi}\,}
\def\dim{\mbox{\rm dim}\,}

\def\dom{\mbox{\rm dom}\,}
\def\ker{\mbox{\rm ker}\,}

\def\cl*co{\mbox{\rm cl}^*\mbox{\rm co}\,}
\def\cl{\mbox{\rm cl}\,}

\def\tto{\rightrightarrows} \def\Hat{\widehat}
\def\h{\hfill\triangle}

\def\Om{\Omega}
\def\O{\Omega}
\def\Th{\Theta}

\def\emp{\emptyset}
\def\st{\stackrel}

\def\Th{\Theta}

\def\N{\mathbb{N}}

\def\h{\hfill\Box}

\begin{document}

\title[Subdifferential Calculus for Ordered Set-Valued Mappings]{\bf Subdifferential Calculus for Ordered Set-Valued Mappings between Infinite-Dimensional Spaces}

\author[B. MORDUKHOVICH]{BORIS S. MORDUKHOVICH}
\thanks{Research of the first author
	was partly supported by the US National Science Foundation under grant DMS-2204519, by the Australian Research Council under Discovery Project DP190100555, and by the Project 111 of China under grant D21024.}
\address[B. MORDUKHOVICH]{\textnormal{Department of Mathematics, Wayne State University, Detroit, MI 48202, USA}}
\email{{\tt aa1086@wayne.edu}} 

\author[O. NGUYEN]{OANH NGUYEN}
\address[O. NGUYEN]{\textnormal{Department of Mathematics and Statistics, Oakland University, Rochester, MI 48309, USA}}
\email{{\tt onguyen@oakland.edu}}

\keywords{Vector and set optimization; variational analysis and generalized differentiation; subdifferentials of ordered multifunctions; subdifferential calculus rules}

\subjclass[2010]{49J52, 49J53, 90C29}

\maketitle
\begin{center}
	\textit{Dedicated to Professor Do Sang Kim, with great respect}
\end{center}

\begin{abstract}
The paper is devoted to developing subdifferential theory for set-valued mappings taking values in ordered infinite-dimensional spaces. This study is motivated by applications to problems of vector and set optimization with various constraints in infinite dimensions. The main results establish new sum and chain rules for major subdifferential constructions associated with ordered set-valued mappings under appropriate qualification and sequentially normal compactness conditions.
\end{abstract}

\section{Introduction}\label{intro}
\setcounter{equation}{0}
\sloppy
This paper belongs to the areas of variational analysis and multiobjective/vector/set optimization in which Professor Do Sang Kim has made many outstanding contributions. We refer the reader to, e.g., the books \cite{KTZ,M1a,M2,TW} and more recent papers \cite{KMPT,kim,tuyen} with the bibliographies therein.

In order to efficiently investigate constrained optimization and equilibrium problems involving set-valued and nonsmooth vectorial mappings with values in ordered spaces, appropriate notions of generalized differentiation for such mappings are needed with the crucial requirement that they satisfy adequate calculus rules. Developing the {\em dual-space approach} to variational analysis and generalized differentiation implemented in the two-volume monograph \cite{M1} for mappings without ordering structures and for extended-real-valued functions, we now focus on extending {\em subdifferential} constructions for set-valued mappings (multifunctions) with values in ordered spaces. Some of such constructions have been introduced and studied in \cite{BM07_1,BM10,M2}, while no calculus rules for them had been obtained therein. In our previous publication \cite{MO}, we developed subdifferential sum and chain rules for the ordered multifunctions between finite-dimensional spaces with applications to the existence of Pareto-type minimizers and deriving necessary optimality conditions. Some calculus results are obtained in \cite{anh,chuong} for subdifferentials of convex mappings associated with efficient point multifunctions in convex vector optimization and convex robust systems in finite-dimensional spaces. Note also related results of \cite{parra,huy} for calculating coderivatives of Pareto front mappings in problems of multiobjective and vector optimization

The main goal of this paper is to derive extended {\em subdifferential calculus rules} for ordered multifunctions between general classes of {\em infinite-dimensional spaces}. Note that calculus rules for infinite-dimensional multifunctions are essentially more involved than in finite dimension. In particular, they require the usage of the so-called {\em sequential normal compactness} conditions and their ordered counterparts, which hold automatically in finite dimensions. The obtained subdifferential calculus rules admit broad applications to the existence theorems and optimality conditions in set-valued optimization with various constraints in infinite-dimensional spaces, while we do not include such results into this paper due to the size limitation.\vspace*{0.05in}

The rest of the paper is organized as follows. In Section~2, we overview some basic tools of variational analysis and generalized differentiation with formulating major facts used in the subsequent material. Section~3 presents the underlying definitions and discussions of subdifferentials of ordered set-valued mappings and related properties that are largely exploited in deriving the main results of  subdifferential calculus for such multifunctions.

In Section~4, we establish fundamental subdifferential sum rules for the four major types of limiting subdifferentials of ordered set-valued mappings under fairly general qualification and sequential normal compactness conditions, which are satisfied for broad classes or ordered multifunctions between infinite-dimensional spaces. The next Section~5 contains several subdifferential chain rules for general compositions of ordered set-valued and single-valued mappings. In the concluding Section~6, we discuss applications of the obtained subdifferential calculus results to problems of set-valued optimization with general constraints of various types.

In this paper, we employ the conventional notation of variational analysis, set-valued analysis, and generalized differentiation; see, e.g., the books \cite{KTZ,M1,M2,RW}.

\section{Machinery of Variational Analysis}\label{prel}
\setcounter{equation}{0}

This sections recalls, following mainly the the monograph \cite{M1}, some basic notions and results of variational analysis and generalized differentiation that are used in the subsequent material. Unless otherwise stated, all the spaces under consideration are {\em normed} equipped with the norm $\|\cdot\|$ and the canonical pairing $\langle \cdot, \cdot \rangle$ between X and its topological dual $X^*$.

Given a set-valued mapping/multifunction $F \colon X \tto X^*$, the (sequential Painlev\'e-Kuratowski) {\em outer limit} of F at $\ox$ in the norm topology of X and weak$^*$ topology $w^*$ of $X^*$ is defined by
\begin{equation}
	 \aligned
\Limsup_{x \to \ox} F(x):= \big\{x ^* \in X^*\;\big|\;\exists\,x_k\to\bar x,\;x_k^* \st{w^*}{\to}x^*\text{ with }x_k^* \in F(x_k)\;\\ \mbox{as}\;k \in \N:=\{1,2,\ldots\}\big\}.
\endaligned\end{equation}

For a nonempty set $\Om\subset X$, the (Fr\'echet) {\em regular normal cone} to $\O$ at $x\in\Om$ is given by
\begin{equation}\label{rnc}
\Hat N(x;\Om):=\Big\{x^*\in X^* \;\Big|\;\limsup_{u\st{\Om}{\to}x}\frac{\langle x^*,u-x\rangle}{\|u-x\|}\le 0\Big\},
\end{equation}
where the symbol `$u\st{\Om}{\to}x$' indicates that $u\to x$ with $u\in\Om$. We let $\Hat N(x;\Om):=\emp$ if $x\notin\Om$. The (limiting, Mordukhovich) {\em normal cone} to $\Om$ at $\ox\in\Om$ is defined by
\begin{equation}\label{nc}
N(\ox;\Om):=\Limsup_{x \to \ox}\Hat{N}(x, \Om).
\end{equation}
We clearly have that $\Hat{N}(\ox;\Om)\subset N(\ox;\Om)$ for any $\Om \subset X$ and $\ox \in \Om$. The equality in this inclusion postulates the property of {\em normal regularity} of $\O$ at $\ox$.

Considering a set-valued mapping $F\colon X\tto Z$ between normed spaces, we associate with it the following sets called the {\em domain}, {\em graph}, and {\em kernel} of $F$, respectively:
\begin{equation*}
\dom F:=\big\{x\in X\;\big|\;F(x)\ne\emp\big\},
\end{equation*}
\begin{equation*}
\gph F:=\big\{(x,z)\in X \times Z\;\big|\;z\in F(x)\big\}.
\end{equation*}
\begin{equation*}
\ker F:=\big\{x \in X\;\big|\;0\in F(x)\big\}.
\end{equation*}
The {\em inverse mapping} of $F \colon X \tto Z$ is $F^{-1}\colon Z \tto X$ with $F^{-1}(y)=\left\{ x \in X \;\big|\; y \in F(x) \right\}.$

Given now a multifunction $F\colon X \tto Z$ and a point $(\ox,\bar z)\in\gph F$, the {\em regular coderivative} of $F$ at $(\ox, \bar z)$ is the set-valued mapping $\Hat{D}^*F(\ox,\oz): Z^* \tto X^*$ with
\begin{equation}\label{reg-cod}
\Hat{D}^*F(\ox,\bar z)(z^*):=\big\{x^*\in X^*\;\big|\;(x^*,-z^*)\in \Hat{N}\big((\ox,\bar z);\gph F\big\},\quad z^*\in Z^*.
\end{equation}
The {\em normal coderivative} of $F$ at $(\ox, \oz)$ is defined by
\begin{equation}\label{2.3}
D_N^*F(\ox,\bar z)(z^*):=\big\{x^*\in X^*\;\big|\;(x^*,-z^*)\in N\big((\ox,\oz);\gph F\big\},\quad z^*\in Z^*,
\end{equation}
and the {\em mixed coderivative} of $F$ at $(\ox, \oz)$ is
\begin{equation}\label{2.4}
\begin{aligned}
D_M^*F(\ox,\oz)(z^*):=&\big\{x^*\in X^*\;\big|\;\exists\,(x_k,z_k) \st{\gph F}{\to} (\ox, \oz),\,x_k^* \st{w^*}{\to} x^*,\,z_k^* \to z^*\\
&\text{ such that }\,(x_k^*,-z_k^*)\in \Hat{N}\big((x_k, z_k);\gph F\big)\big\},\quad z^*\in Z^*.
\end{aligned}
\end{equation}
If $F(\ox)$ is a singleton at $\ox$, we skip $\oz=F(\ox)$ in the coderivative notations.

The difference between the mixed coderivative and normal coderivatives is that in \eqref{2.3} we use the weak$^*$ convergence for both sequences $\{x_k^*\}$ and $\{z_k^*\}$, while \eqref{2.4} employs the weak$^*$ convergence for $\{x_k^*\}$ and the norm convergence for $\{z_k^*\}$ as $k\to\infty$.

We obviously  have the inclusions
\begin{equation}\label{2.5}
\Hat{D}^*F(\ox,\oz)(z^*) \subset D_M^*F(\ox,\oz)(z^*) \subset D_N^*F(\ox,\oz)(z^*) \;\text{ for all }\;z^* \in Z^*,
\end{equation}
where the rightmost one holds as an equality if $Z$ is finite-dimensional.

Recall that a single-valued mapping $f: X \to Z$ is {\em strictly differentiable} at $\ox$ with the strict derivative $\nabla f(\ox)$ of $f$ at this point if
\begin{equation*}
\lim_{x,u\to\ox}\frac{f(x)-f(u)-\langle\nabla f(\ox),x-u\rangle}{\|x-u\|}=0,
\end{equation*}
which happens, in particular, when $f$ is continuously differentiable $(\mathcal{C}^1$-smooth) around the point in question. If $f$ is strictly differentiable at $\ox$, then
\begin{equation}
D_M^*F(\ox,\oz)(z^*) = D_N^*F(\ox,\oz)(z^*)= \big\{ \nabla f(\ox)^*z^*\big\}\;\text{ for all }\;z^* \in Z^*.
\end{equation}

Similarly to the normal regularity of sets, we define the {\em graphical regularity} of mappings. To be more specific, it is said that a set-valued mapping $F$ is {\em N-regular} at $(\ox, \oz) \in \gph F$ if $D_N^*F(\ox, \oz)=\Hat{D}^*F(\ox, \oz)$.  It is {\em M-regular} at this point if $D_M^*F(\ox, \oz)=\Hat{D}^*F(\ox, \oz)$. We obviously have that the normal regularity of $F$ implies its mixed regularity at this point but not vice versa.

Next we recall fundamental {\em well-posedness} properties of set-valued mappings, which play a crucial role in nonlinear and variational analysis as well as in their numerous applications. A multifunction $F\colon X\tto Z$ with $\dom F \neq \emp$ is {\em Lipschitz-like} around $(\ox,\oz)\in\gph F$ with modulus $\ell\ge 0$ if there exist neighborhoods $U\subset X$ of $\ox$ and $V \subset Z$ of $\oz$ such that
\begin{equation}\label{lip-like}
F(x)\cap V\subset F(u)+\ell\,\|x-u\|\mathbb B\;\mbox{ for all }\;x,u\in U,
\end{equation}
where $\mathbb B$ stands for the closed unit ball in the space in question. If $F=f\colon X\to Z$ is single-valued, then \eqref{lip-like} reduces go the classical local Lipschitzian property of $f$ around $\ox$:
\begin{equation*}
\|f(x)-f(u)\| \le\ell\|x-u\|\;\mbox{ for all }\;x,u\in  U.
\end{equation*}
It is easy to see that any mapping $f\colon X\to Z$, which is strictly differentiable at $\ox$, is Lipschitz continuous around this point.

The other fundamental property  of set-valued mappings used below is defined as follows. A multifunction $F\colon X\tto Z$ is {\em metrically regular} around $(\ox,\oz)\in\gph F$ with modulus $\mu>0$ if there exist neighborhoods $U \subset X$ of $\ox$ and $V$ of $\oz$ such that
\begin{equation}\label{mr}
{\rm dist}\big(x;F^{-1}(z)\big)\le\mu\,{\rm dist}\big(z;F(x)\big)\;\mbox{ for all }\;x\in U\;\mbox{ and }\;z\in V,
\end{equation}
where `dist' signifies the distance between a point and a set.

It has been well recognized in variational analysis (see, e.g., \cite[Theorem~1.49]{M1}) that $F\colon X\tto Z$ is Lipschitz-like around $(\ox, \oz)$ {\em if and only if} the inverse mapping $F^{-1}\colon Z \tto X$ is metrically regular around $(\oz, \ox)\in \gph F^{-1}$ with the same exact bounds of the corresponding moduli.\vspace*{0.05in}

The main results of this paper are obtained in the case of Asplund spaces. It is said that a Banach space $X$ is {\em Asplund} of each of its separable subspace has a separable dual. This subclass of Banach spaces is sufficiently broad; in particular, it contains every reflexive Banach space; see, e.g., \cite{fabian,M1} for more facts, references, and discussions. We refer the reader to the books \cite{M1,M1a,M2,M24} and the bibliographies therein for a variety of results on generalized differentiation in Asplund spaces and their applications to different aspects of variational analysis, optimization, and control.

Before formulating the major characterizations of the Lipschitz-like and hence metric regularity properties of multifunctions between Asplund spaces, we recall that a certain amount of ``normal compactness'' is required to furnish limiting procedures, which involve weak$^*$ sequential limits of normals in infinite-dimensional spaces.

A nonempty subset $\O\subset X$ of a normed space is {\em sequentially normally compact} at $\ox\in\O$ if for any sequences $x_k\st{\O}{\to}\ox$ and $x^*_k\st{w^*}{\to}0$ with $x^*_k\in\Hat N(x_k;\O)$  for all $k\in\N$ we have $\|x^*_k\|\to 0$ as $k\to\infty$. This property is obviously automatic in finite-dimensional spaces. On the other hand, it has been recognized in variational analysis that the SNC property holds for broad classes of sets in infinite dimensions; in particular, those which possess a certain Lipschitzian behavior; see \cite{M1}.

The weakest property of this type used below for the case of multifunctions, is defined as follows. A set-valued mapping $F\colon X\tto Z$ is {\em partially sequentially normally compact} (PSNC) at $(\ox,\oz)\in\gph F$ if for any $(x_k, z_k)\in X\times Z$, $(x_k^*, z_k^*) \in X^* \times Z^*$ such that $(x_k, z_k)\st{\gph F}{\to} (\ox, \oz)$ and $x_k^* \in \Hat{D}^*F(x_k,z_k)(z_k^*)$ whenever $k \in \N$, we have the implication
\begin{equation}\label{psnc}
\big[x_k^*\st{w^*}{\to} 0,\;\|z_k^*\|\to 0\big] \implies \|x_k^*\|\to 0\;\text{ as }\;k \to \infty.
\end{equation}
This property obviously holds if $X$ is a finite-dimensional space. It follows from Theorem~\ref{lipschitz like} below that any closed-graph mapping that is Lipschitz-like around $(\ox,\oz)\in\gph F$, is PSNC at this point. More results and applications for this property, including efficient conditions for its preservation (calculus) under various operations on multifunctions, can be found in \cite{M1}.\vspace*{0.05in}

Now we are ready to formulate the major characterization of the Lipschitz-like property of set-valued mappings between Asplund spaces in terms of the mixed coderivative \eqref{2.4} and the the PSNC property of multifunctions \eqref{psnc}.

\begin{theorem}[\bf characterization of Lipschitz-like multifunctions]\label{lipschitz like} Let $F \colon X \tto Z$ be a set-valued mapping between Asplund spaces, and let the graph of $F$ be closed around $(\ox, \oz) \in \gph F$. Then the following properties are equivalent:
\begin{enumerate}
\item[\bf(i)] $F$ is Lipschitz-like around $(\ox, \oz)$.
\item[\bf(ii)] $F$ is PSNC at $(\ox, \oz)$ and $D^*_M F(\ox,\oz)(0)=\left\{0 \right\}$.
\end{enumerate}
\end{theorem}

In finite dimensions, the PSNC property holds automatically, and the coderivative characterization of Theorem~\ref{lipschitz like} reduces to the result of \cite{M93} named in \cite{RW} the {\em Mordukhovich criterion}. Observe that in the general setting of Asplund spaces, the coderivative and PSNC constructions of Theorem~\ref{lipschitz like} enjoy {\em full calculi} developed in \cite{M1}, which ensure broad applications of this theorem for large classes if structured problems. Note also that, by the aforementioned equivalence between the Lipschitz-like and metric regularity properties of multifunctions, Theorem~\ref{lipschitz like} yields the corresponding characterizations of metric regularity; see  \cite[Theorem~4.18]{M1}.

\section{Ordered Set-Valued Mappings and Their Subdifferentials}
\setcounter{equation}{0}

In this section, we consider multifunctions with values in {\em ordered} normed spaces and define some {\em subdifferential} constructions for them. Given a normed space $Z$ and a nonempty, closed, convex cone $\Th\subset Z$, define the {\em preference relation} $\preceq$ on $Z$ by
\begin{equation}\label{order}
z_1\preceq z_2\Longleftrightarrow z_2-z_1\in\Theta.
\end{equation}
For a set-valued mapping $F\colon X\tto Z$ on a normed space $X$ with values in a normed space $Z$ ordered by \eqref{order}, the {\em generalized epigraph} of $F$ with respect to $\Theta$ is defined by
\begin{equation}\label{epi}
\epi_\Theta F\:=\big\{(x,z)\in X\times Z\;\big|\;z\in F(x)+\Theta\big\}
\end{equation}
and the {\em $\Th$-epigraphical multifunction} ${\mathcal{E}}_{F,\Theta}\colon X\tto Z$ is given by
\begin{equation}\label{epim}
{\mathcal{E}}_{F,\Theta}(x):=\big\{z\in Z\;\big|\;z\in F(x)+\Theta\big\}.
\end{equation}
It is easy to see that $\gph{\mathcal{E}}_{F,\Theta}=\epi_\Theta F$.  We say that  a multifunction $F\colon X\tto Z$  with $\Th$-ordered values is {\em $\Th$-epiclosed} if its generalized epigraph \eqref{epi} is closed in $X\times Z$.

Induced by the notions of graphical regularity of arbitrary multifunctions introduced in Section~2, we say that a set-valued mapping $F\colon X\tto Z$ with ordered values in $Z$ is is {\em $N$-epiregular} (resp.\ {\em $M$-epiregular}) at $(\ox, \oz) \in \epi_\Theta F$ if the corresponding $\Th$-epigraphical multifunction ${\mathcal{E}}_{F,\Theta}$ is $N$-regular (resp.\ $M$-regular) at this point.\vspace*{0.05in}

Next  we introduce the subdifferential notions for ordered multifunctions, which reduce to the corresponding subdifferentials of scalar (extended-real-valued) functions  when \eqref{order} signifies the standard order on the real line; cf.\ \cite{M1}.

\begin{definition}[\bf subdifferentials of ordered set-valued mappings]\label{sub_or} Let $F\colon X\tto Z$ be a multifunction between normed spaces, where $Z$ is ordered by the preference in \eqref{order} generated via a closed and convex cone $\Th\ne\emp$. Having $(\ox,\oz)\in\epi_\Theta F$ from \eqref{epi} and using the coderivatives in \eqref{reg-cod}--\eqref{2.4} for the epigraphical multifunction \eqref{epim}, define the following constructions:
\begin{itemize}
\item[]{\bf(i)} The {\sc regular subdifferential}  of $F$ at $(\ox,\oz)\in \epi_{\Theta}F$ is
\begin{equation*}
\Hat{{\partial}}_{\Theta} F(\ox,\oz):=\big\{x^*\in X^*\;\big|\;x^*\in\Hat{D}^*{\mathcal{E}}_{F,\Theta}(\ox,\oz)(z^*),\;-z^*\in N(0;\Theta),\;\|z^*\|=1\big\}.
\end{equation*}
\item[]{\bf(ii)} The {\sc normal subdifferential} of $F$ at $(\ox,\oz)\in \epi_{\Theta}F$ is
\begin{equation*}
{\partial}_{N,\Theta} F(\ox,\oz):=\big\{x^*\in X^*\;\big|\;x^*\in{D}_N^*{\mathcal{E}}_{F,\Theta}(\ox,\oz)(z^*),\;-z^*\in N(0;\Theta),\;\|z^*\|=1\big\}.
\end{equation*}
\item[]{\bf(iii)} The {\sc mixed subdifferential} of $F$ at $(\ox,\oz)\in \epi_{\Theta}F$ is
\begin{equation*}
{\partial}_{M,\Theta} F(\ox,\oz):=\big\{x^*\in X^*\;\big|\;x^*\in{D}_M^*{\mathcal{E}}_{F,\Theta}(\ox,\oz)(z^*),\;-z^*\in N(0;\Theta),\;\|z^*\|=1\big\}.
\end{equation*}
\item[]{\bf(iv)} The {\sc normal singular subdifferential} of F at $(\ox,\oz)\in \epi_{\Theta}F$ is
\begin{equation*}
{\partial}^\infty_{N,\Theta} F(\ox,\oz):={D}_N^*{\mathcal{E}}_{F,\Theta}(\ox,\oz)(0).
\end{equation*}
\item[]{\bf(v)} The {\sc mixed singular subdifferential} of F at $(\ox,\oz)\in \epi_{\Theta}F$ is
\begin{equation*}
{\partial}^\infty_{M,\Theta} F(\ox,\oz):={D}_M^*{\mathcal{E}}_{F,\Theta}(\ox,\oz)(0).
\end{equation*}
\end{itemize}
\end{definition}

It the epigraphical multifunction \eqref{epim} is $N$-regular as defined in Section~2 (in particular, when $\dim Z<\infty$), then it follows from the equality in the rightmost inclusion of \eqref{2.5} that the normal subdifferential and the mixed subdifferential agree. The same identity hold for the normal and mixed singular subdifferentials of ordered multifunctions from Definition~\ref{sub_or}.

To proceed further, let us formulate a version of the PSNC property for ordered multifunctions. Given $F\colon X\tto Z$ having values in an ordered space $Z$ with the ordering cone $\Th$  in \eqref{order}, we say that $F$ is {\em partially sequentially normally epi-compact} (PSNEC) at $(\ox, \oz)\in \epi_\Theta F$ if the epigraphical multifunction ${\mathcal{E}}_{F,\Theta}$ from \eqref{epim} is PSNC at this point.

Finally in this section, we say that an ordered mapping $F\colon X\tto Z$ is {\em epi-Lipschitz-like} (ELL) around $(\ox, \oz)$ with modulus $\ell\ge 0$ if the epigraphical multifunction  ${\mathcal{E}}_{F,\Theta}$ is Lipschitz-like around this point with modulus $\ell\ge 0$, i.e., if there exist neighborhoods $U$ of $\ox$ and $V$ of $\oz$ such that
\begin{equation}\label{ell}
{\mathcal{E}}_{F, \Theta} (x) \cap V \subset {\mathcal{E}}_{F,\Theta}(u)+\ell ||x-u||\mathbb B\; \mbox{ for all }\;x, u \in U.
\end{equation}
It follows from Theorem~\ref{lipschitz like} and Definition~\ref{sub_or}(v) that for ordered multifunctions $F\colon X\tto Z$ with closed generalized epigraph \eqref{epi} that the simultaneous fulfillment of the PSNEC property of $F$ at $(\ox,\oz)\in \epi_\Theta F$ and the mixed singular subdifferential condition
\begin{equation}\label{ell-singular}
{\partial}^\infty_{M,\Theta} F(\ox,\oz)=\{0\}
\end{equation}
provides a {\em complete characterization} of the ELL property \eqref{ell} of $F$ around $(\ox,\oz)$.

\section{Sudifferential Sum Rules for Ordered Set-Valued Mappings}
\setcounter{equation}{0}

This section is mainly devoted to deriving fundamental {\em sum riles} for the normal and mixed subdifferentials, as well as for their singular counterparts from Definition~\ref{sub_or}, of ordered multifunctions between Asplund spaces.

We begin with recalling the following properties for arbitrary multifunctions $S\colon X\tto Y$ between normed spaces. Given a point $(\ox,\bar y)\in\gph S$, the set-valued mapping $S$ is said to be {\em inner semicontinuous} at $(\ox,\bar y)$ if
\begin{equation}\label{inner}
\Liminf_{x \to \ox} S(x)= \bar y \in S(\ox),
\end{equation}
i.e., by using the {\em inner limit} construction ``Liminf" taken from, e.g., \cite{M1,RW},  for every sequence $x_k\to\ox$, there exists a sequence $z_k\in S(x_k)$ that converges to $\bar y$ as $k\to\infty$. The multifunction $S$ is
{\em inner semicompact} at $\ox\in\dom S$ if for every sequence $x_k \to \ox$, there is a sequence $y_k \in S(x_k)$ that contains a convergent subsequence as $k \to \infty$.

The inner semicontinuity property \eqref{inner} obviously holds for multifunctions that are Lipschitz-like around $(\ox,\bar y)$. Indeed, the latter property yields even the continuity of $S$ around $(\ox,\bar y)$, i.e., that the outer and inner limits of $S$ agree for all points near the reference one $(\ox,\bar y)$. The inner semicompactness is significantly more general than the inner semicontinuity of the multifunction while not requiring any amount of (semi)continuity. In particular, for multifunctions with values in finite-dimensional spaces, the {\em uniform boundedness} of $S(x)$ around $\ox$ yields the inner semicompactness of $S$ at the point in question.\vspace*{0.05in}

Given now two set-valued mappings $F_1, F_2 \colon X\tto Z$ between normed spaces, define the  multifunction $S \colon X\times Z \tto Z\times Z$ associated with $F_1,F_2$ by
\begin{equation}\label{3.1}
S(x,z):= \big\{(z_1,z_2) \in Z\times Z\;\big|\;z_1 \in F_1(x),\;z_2 \in F_2(x),\; z=z_1+z_2\big\}.
\end{equation}

Before establishing the main result of this section, we present two lemmas. The first one, taken from  \cite[Theorem~3.10]{M1}, gives us some versions of the coderivative sum rules for set-valued mappings between Asplund spaces. The presented versions depend on whether mapping \eqref{3.1} is assumed to be lower semicontinuous or lower semicompactness at the point in question. Observe that the coderivative sum rules are obtained in parallel for both normal \eqref{2.3} and mixed \eqref{2.4} coderivatives under the qualification condition expressed in both cases via the mixed coderivative.

\begin{lemma}[\bf coderivative sum rules]\label{coder_sum} Let $F_1, F_2 \colon X\tto Z$ be sert-valued mappings between Asplund spaces,  and let $(\ox,\oz)\in\gph(F_1+F_2)$. The following two independent assertions hold:
\begin{itemize}
\item[] {\bf (i)} Fix $(\oz_1,\oz_2) \in S(\ox,\oz)$ from \eqref{3.1} and suppose that this mapping is inner semicontinuous at $(\ox,\oz,\oz_1,\oz_2)$. Assume also that the graphs of $F_1$ and $F_2$ are locally closed around $(\ox, \oz_1)$ and $(\ox, \oz_2)$, respectively, that either $F_1$ is PSNC at  $(\ox, \oz_1)$ or $F_2$ is PSNC at  $(\ox, \oz_2)$, and that the mixed coderivative qualification condition
\begin{equation}\label{qc_sum}
D_M^*F_1(\ox,\oz_1)(0)\cap\big(-D_M^*F_2(\ox,\oz_2)(0)\big)=\{0\}
\end{equation}
is satisfied. Then for all $z^*\in Z^*$ we have the inclusion
\begin{equation}\label{coder_sum_subset}
D^*(F_1+F_2)(\ox,\oz)(z^*)\subset D^*F_1(\ox,\oz_1)(z^*)+D^*F_2(\ox,\oz_2)(z^*),
\end{equation}
where $D^*$ in both parts stands for either normal or mixed coderivative.
		
\item [] {\bf(ii)} If instead of the inner semicontinuity of \eqref{3.1} at $(\ox,\oz)\in\gph S$, we assume that $S$ is inner semicompact at $\ox\in\dom S$ and that all the conditions of ${\rm(i)}$ are satisfied for every $(\oz_1,\oz_2)\in S(\ox,\oz)$, then whenever $z^*\in Z^*$ we have the inclusion
\begin{equation*}
D^*(F_1+F_2)(\ox,\oz)(z^*)\subset \bigcup\limits_{(\oz_1,\oz_2)\in S(\ox,\oz)}\big[D^*F_1(\ox,\oz_1)(z^*)+D^*F_2(\ox,\oz_2)(z^*)\big]
\end{equation*}
for both normal and mixed coderivatives $D^*=D^*_N,D^*_M$.
\end{itemize}
\end{lemma}

To proceed further with developing {\em subdifferential sum rules} for {\em ordered} multifunctions, we need the next lemma providing the {\em epigraphical sum rule} for such multifunctions between arbitrary normed spaces.

\begin{lemma}[\bf epigraphical multifunctions under set summation]\label{epi_sumrule} Let $F_1,F_2\colon X\tto Z$ be set-valued mappings between normed spaces, where $Z$ is ordered by a closed convex cone $\Th\subset Z$ in \eqref{order}. Then for any $x\in X $, we have the equality
\begin{equation}
{\mathcal{E}}_{F_1+F_2,\Theta}(x)={\mathcal{E}}_{F_1,\Theta}(x)+{\mathcal{E}}_{F_2,\Theta}(x),
\end{equation}
where ${\mathcal{E}}$ is the epigraphical multifunction defined in \eqref{epim}.
\end{lemma}
{\noindent\it Proof}. Take any $x \in X$ and $z\in{\mathcal{E}}_{F_1+F_2,\Theta}(x)=(F_1+F_2)(x)+\Theta$. Then there exists $z_1\in F_1(x)+\Theta$ and $z_2\in F_2(x)+\Theta$ such that $z=z_1+z_2$. Therefore, ${\mathcal{E}}_{F_1+F_2,\Theta}(x)\subset {\mathcal{E}}_{F_1,\Theta}(x)+{\mathcal{E}}_{F_2,\Theta}(x)$.

To verify the reverse implication, fix $x\in X$ and take $z \in {\mathcal{E}}_{F_1,\Theta}(x)+{\mathcal{E}}_{F_2,\Theta}(x)$. This means that there exist $z_1\in {\mathcal{E}}_{F_1,\Theta}(x)$ and $z_2\in{\mathcal{E}}_{F_2,\Theta}(x)$ such that $z=z_1+z_2$. Then we find $z_1\in F_1(x)+\Theta$ and $z_2\in F_2(x)+\Theta$ directly by definition \eqref{epim}. This tells us that $z_1+z_2\in {\mathcal{E}}_{F_1+F_2,\Theta}(x)$, which therefore completes the proof of the lemma. $\h$\vspace*{0.1in}

Now we are ready to establish the desired sum rules for the normal and mixed subdifferentials of multifunctions with values in ordered Asplund spaces. Given two set-valued mappings $F_1,F_2\colon X\tto Z$ of this type, consider the mapping $S\colon X\times Z\tto Z\times Z$ defined by
\begin{equation}\label{new_mapping_S}
S_{\mathcal{E}}(x,z):=\big\{(z_1,z_2)\in Z\times Z\;\big|\;z_1\in {\mathcal{E}}_{F_1,\Theta}(x),\;z_2\in {\mathcal{E}}_{F_2,\Theta}(x),\; z=z_1+z_2\big\}
\end{equation}
via the epigraphical multifunctions \eqref{epim} of $F_1$ and $F_2$.\vspace*{0.05in}

Based on the above lemmas and the subdifferential definitions for ordered set-valued mappings, we present {\em eight versions} of the subdifferential sum rules for both normal and mixed subdifferentials from Definition~\ref{sub_or}(ii,iii) and its singular counterparts from Definition~\ref{sub_or}(iv.v) under the semicontinuity and semicompactness assumptions imposed on the multifunction $S_{\mathcal E}$ from \eqref{new_mapping_S}. Let us emphasize here that all the eight subdifferential sum rules are obtained under the same qualification conditions expressed in terms of the {\em mixed singular} subdifferential.

\begin{theorem}[\bf sum rules for subdifferentials]\label{ordered_sum} Let $F_1,F_2:X\tto Z$ be set-valued mappings between Asplund spaces, where $Z$ is ordered by the ordering convex and closed cone $\Th$ in \eqref{order}. Given $(\ox,\oz)\in\epi_\Th(F_1+F_2)$, we get the following assertions:
\begin{itemize}
\item [] {\bf(i)} Fix $(\oz_1,\oz_2)\in S_{\mathcal{E}}(\ox,\oz)$ from \eqref{new_mapping_S} and suppose that $S_{\mathcal{E}}$ is inner semicontinuous at $(\ox,\oz,\oz_1,\oz_2)$.  Assume that $F_1, F_2$ are epiclosed around the corresponding points, and that either $F_1$ is PSNEC at $(\ox, \oz_1)$ or $F_2$ is PSNEC at $(\ox, \oz_2)$. Impose also the mixed singular subdifferential qualification condition
\begin{equation}\label{qc_order_sum}
{\partial}^\infty_{M,\Theta} F_1(\ox,\oz_1)\cap\big(-{\partial}^\infty_{M,\Theta} F_2(\ox,\oz_2)\big)=\{0\}.
\end{equation}
Then we have the subdifferential sum rule inclusions
\begin{equation}\label{order_sum_subset}
\partial_\Theta(F_1+F_2)(\ox,\oz)\subset\partial_\Theta F_1(\ox,\oz_1)+\partial_\Theta F_2(\ox,\oz_2),
\end{equation}
\begin{equation}\label{order_sum_sing}
\partial^\infty_\Theta(F_1+F_2)(\ox,\oz)\subset\partial^\infty_\Theta F_1(\ox,\oz_1)+\partial^\infty_\Theta F_2(\ox,\oz_2)
\end{equation}
valid for both normal and mixed subdifferentials $\partial_\Th=\partial_{N,\Th},\; \partial_{\Th,M}$ in \eqref{order_sum_subset} and the corresponding normal and mixed singular subdifferentials $\partial_\Th^\infty=\partial_{N,\Th}^\infty,\;\partial_{\Th,M}^\infty$ in \eqref{order_sum_sing}.

\item[] {\bf(ii)} If instead of the inner semicontinuity of $S_{\mathcal E}$ in ${\rm(i)} $, we assume that this mapping is inner semicompact at $(\ox,\oz)\in\dom S_{\mathcal E}$ and all the other conditions in ${\rm(i)}$  hold for every $(\oz_1,\oz_2)\in S_{\mathcal E}(\ox,\oz)$, then we have extended  subdifferential sum rule inclusions
\begin{equation*}
{\partial}_\Theta (F_1+F_2)(\ox,\oz)\subset\bigcup\limits_{(\oz_1,\oz_2)\in S_{\mathcal{E}}(\ox,\oz)}\big[{\partial}_\Theta F_1(\ox,\oz_1)+{\partial}_\Theta F_2(\ox,\oz_2)\big],
\end{equation*}
\begin{equation*}
\partial^\infty_\Theta (F_1+F_2)(\ox,\oz)\subset\bigcup\limits_{(\oz_1,\oz_2)\in S_{\mathcal{E}}(\ox,\oz)}\big[\partial^\infty_\Theta F_1(\ox,\oz_1)+\partial^\infty_\Theta F_2(\ox,\oz_2)\big]
\end{equation*}
valid, respectively, for the four subdifferentials discussed above.
\end{itemize}
\end{theorem}
{\noindent\it Proof}. To verify assertion (i), observe first that if $F_1$ and $F_2$ are PSNEC at $(\ox, \oz_1)$ and $(\ox, \oz_2)$, respectively, then ${\mathcal{E}}_{F_1, \Theta}$ and ${\mathcal{E}}_{F_2, \Theta}$ are PSNC at the corresponding points. It is easy to see that the sets $\gph{\mathcal{E}}_{F_1, \Theta}$ and $\gph{\mathcal{E}}_{F_2, \Theta}$ are closed by the imposed epiclosedness of the mappings $F_1$ and $ F_2$. Now we employ the coderivative sum rules from Lemma~\ref{coder_sum}(i) for both normal and mixed coderivatives of the epigraphical multifunctions. Observe that the mixed subdifferential condition \eqref{qc_order_sum} yields the fulfillment of the mixed coderivative qualification condition \eqref{qc_sum} for the epigraphical multifunctions.
Picking any $x^*\in{\partial}_\Theta (F_1+F_2)(\ox,\oz)$ for either
$\partial_\Theta=\partial_{N,\Theta}$ or $\partial_\Theta=\partial_{M,\Theta}$ and then using Lemma~\ref{coder_sum}(i) and Lemma~\ref{epi_sumrule} give us
$$
x^*\in{D}^*{\mathcal{E}}_{F_1+F_2,\Theta}(\ox,\oz)(z^*)={D}^*({\mathcal{E}}_{F_1,\Theta}+{\mathcal{E}}_{F_2,\Theta})(\ox,\oz)(z^*)
$$
with $-z^*\in N(0;\Theta)$ and $\|z^*\|=1$. Applying the sum rule for coderivatives, under the inner semicontinuity property of $S_{\mathcal{E}}$ and the qualification condition \eqref{qc_sum} for the epigraphical multifunctions, brings us to the inclusion
$$
x^*\in {D}^*{\mathcal{E}}_{F_1,\Theta}(\ox,\oz_1)(z^*)+D^*{\mathcal{E}}_{F_2,\Theta}(\ox,\oz_2)(z^*),
$$
which yields $x^* \in \partial_\Theta F_1(\ox, \oz_1)+\partial_\Theta F_2(\ox,\oz_2)$ and thus justifies assertions (i).

To verify  now assertion (ii) of the theorem, we apply the coderivative sum rules from  Lemma~\ref{coder_sum}(ii) to the epigraphical multifunctions associated with $F_1$ and $F_2$. Proceeding similarly  to the proof of (i) with the usage of the inner compactness of $S_{\mathcal{E}}$ at $(\ox, \oz)$ together with the imposed qualification condition and the PSNC properties of the mappings ${\mathcal{E}}_{F_1,\Theta}, {\mathcal{E}}_{F_2,\Theta}$ at the corresponding points, we arrive at the subdifferential sum rules in (ii) and thus complete the proof of the theorem. $\h$\vspace*{0.05in}

Note that the mixed subdifferential qualification condition \eqref{qc_order_sum} can be replaced in Theorem~\ref{ordered_sum} by the normal subdifferential qualification condition
\begin{equation}\label{normal-qc}
\partial^\infty_{N,\Theta}F_1(\ox,\oz_1)\cap\big(-\partial^\infty_{N,\Theta}F_2(\ox,\oz_2)\big)=\{0\},
\end{equation}
which is however more restrictive in infinite dimensions. In particular, the following corollary is a consequence of Theorem~\ref{ordered_sum}(i) under the mixed subdifferential qualification condition \eqref{qc_sum} but not under its normal counterpart \eqref{normal-qc}.

\begin{corollary}[\bf subdifferential sum rules for ELL mappings]\label{sum-lip} Let $F_1,F_2\colon X\tto Z$ be set-valued mappings between Asplund spaces, where $Z$ is ordered by closed and convex cone $\Th$ in \eqref{order}, and let $(\ox,\oz) \in \epi_\Theta (F_1+F_2)$. Fix $(\oz_1,\oz_2)\in S_{\mathcal{E}}(\ox, \oz)$ such that $S_{\mathcal{E}}$ is inner semicontinuous at $(\ox, \oz,\oz_1,\oz_2)$ and suppose that $F_1$ and $F_2$ are epiclosed around $(\ox, \oz_1)$ and $(\ox,\oz_2)$, respectively. Assume further that either $F_1$ is ELL around $(\ox, \oz_1)$ or $F_2$ is ELL around $(\ox,\oz_2)$. Then the sum rules  in \eqref{order_sum_subset} hold for both normal and mixed subdifferentials $\partial_\Th=\partial_{N,\Th},\;\partial_{M,\Th}$, while those in \eqref{order_sum_sing} are fulfilled for their singular counterparts $\partial_\Th^\infty=\partial^\infty_{N,\Th},\;\partial^\infty_{M,\Th}$.
\end{corollary}
{\noindent\it Proof}. Suppose for definiteness that $F_1$ is ELL around $(\ox, \oz_1)\in \epi_\Theta F_1$. Then it follows from \eqref{ell-singular} that ${\partial}^\infty_{M,\Theta} F_1 (\ox,\oz)=\{0\}$, which justifies the qualification condition \eqref{qc_order_sum}. On the other hand, we have that the epigraphical multifunction ${\mathcal{E}}_{F_1,\Th}$ is locally Lipschitz-like around $(\ox, \oz_1)\in \gph {\mathcal{E}}_{F_1, \Th}=\epi_\Theta F_1$, and therefore ${\mathcal{E}}_{F_1,\Th}$ is PSNC at $(\ox, \oz_1)$ by Theorem~\ref{lipschitz like}. This shows that all the requirements of Theorem~\ref{ordered_sum}(i) hold, and thus we are done with the proof. $\h$\vspace*{0.05in}

Given now a nonempty set $\O\subset X$ in a normed space $X$, we associate with it the {\em indicator mapping} $\Delta(\cdot; \Om)\colon X \tto Z$ relative to another normed space $Z$ by
\begin{equation}\label{ind-map}
\Delta(x;\Omega):=\left\{\begin{array}{ll}
0&\mbox{if }\;x\in\Omega,\\
\emp&\mbox{if }\;x\notin\Omega.
\end{array}\right.
\end{equation}
Observe the useful subdifferential relationships for the indicator mapping \eqref{ind-map} given below.

\begin{lemma}[\bf singular subdifferentials of the indicator mappings]\label{sub indicator} Let $\Om\subset X$ be a subset of a normed space with $\ox\in\O$, and let $\Delta(\cdot;\Om)\colon X \tto Z$ be the indicator mapping of $\O$ relative to a normed space $Z$ ordered by a closed convex cone $\Th\subset Z$ from \eqref{order}. Then we have the following relationships for the normal and mixed singular subdifferential of $\Delta(\cdot;\Om)$ at $(\ox,0)\in\gph\Delta(\cdot;\Om)$ and the normal cone to  the generating set $\O$:
\begin{equation}\label{3.13}
\partial_{M,\Theta}^\infty \Delta(\cdot;\Om)(\ox; 0)=\partial_{N,\Theta}^\infty \Delta(\cdot;\Om)(\ox, 0)= N(\ox; \Om).
\end{equation}
\end{lemma}
{\noindent\it Proof}. Fix $\ox \in \Om$. The first equality in \eqref{3.13} follows directly from the definitions of the normal and mixed subdifferentials via the coderivatives of the indicator mapping \eqref{ind-map}. To verify the second equality in \eqref{3.13}, take an arbitrary element $x^*\in \partial_{N,\Theta}^\infty\Delta(\cdot; \Om)(\ox, 0)$ and deduce from Definition~\ref{sub_or}(iv) together with the one for the normal coderivative \eqref{2.3} that we equivalently have $(x^*,0)\in N((\ox, 0);\epi_\Theta \Delta(\cdot,\Om))$. Note that
$$
\epi_\Theta \Delta(\cdot; \Om)=\big\{(x,z)\;\big|\;z \in \Delta(x;\Om)+\Theta\big\}=\Om\times\Theta,
$$
which shows that $(x^*,0)\in N((\ox,0);\Om\times\Theta)= N(\ox;\Om)\times N(0;\Theta)$ by the product rule for normal cones valid in any normed space; see, e.g., \cite[Proposition~1.2]{M1}. Therefore, $x^* \in N(\ox;\Om)$, which completes the proof of the lemma. $\h$\vspace*{0.05in}

The next lemma provides a relationship between the regular subdifferential of an ordered multifunction from Definition~\ref{sub_or}(i) and its restriction on a given set.

\begin{lemma}[\bf regular subdifferentials of the restricted multifunctions]\label{reg sub f and indicator} Let $F\colon X\tto Z$ be a set-valued mapping between normed spaces, where $Z$ is ordered by the closed convex cone $\Th\subset Z$ in \eqref{order}. Define the restricted multifunction $F_\O\colon X\tto Z$ by
\begin{equation}\label{restr}
F_\O(x):=F(x)+\Delta(x;\O),\quad x\in\O,
\end{equation}
via the indicator mapping \eqref{ind-map} associated with a given set $\O$. Then we have the inclusion
\begin{equation}
\Hat{\partial}_\Theta F(\ox, \oz)+\Hat {N}(\ox, \Om)\subset \Hat{\partial}_\Theta F_\Om (\ox,\oz).
\end{equation}
\end{lemma}
{\noindent\it Proof}. Observe first that $\epi_\Theta F_\Om=\epi_\Theta F \cap (\Om \times Z)$. Indeed, picking $(x, z) \in \epi_\Theta F_\Om$ means that $z \in F(x)+\Theta$, where $x \in \Om$. This is clearly equivalent to $(x,z) \in \epi_\Theta F\cap(\Om\times Z)$.

Further, take $x^*=u^*+v^* \in\Hat{\partial}_\Theta F(\ox, \oz)+\Hat{N}(\ox, \Om)$. where $u^* \in \Hat{\partial}_\Theta F(\ox, \oz)$ and $v^*\in\Hat{N}(\ox, \Om)$. It follows from the regular subdifferential construction and the aforementioned product rule for the  normal cones involved that
$$
(u^*, -z^* )\in \Hat{N}((\ox, \oz);\epi_\Theta F)\;\text{ with }\;-z^*\in N(0;\Theta)\;\mbox{ and }\;\|z^*\| =1$$
together with the regular normal cone inclusion
$$
(v^*, 0)\in \Hat{N}\big((\ox, \oz);\Om \times Z\big).
$$
Therefore, by adding to both sides the above relationship between the generalized epigraphs $\epi_\Theta F_\Om$ and $\epi_\Theta F$, we arrive at the inclusion
$$
(u^*+v^*, -z^*) \in \Hat{N}\big((\ox, \oz); \epi_\Theta F\big)+\Hat{N}\big((\ox, \oz);\Om\times Z\big)\subset \Hat{N}\big((\ox, \oz);\epi_\Theta F_\Om\big).
$$
Combining the latter with the conditions imposed on $z^*$ tells us that $x^*=u^*+v^* \in \Hat{\partial}_\Theta F_\Om(\ox, \oz)$, which thus completes the proof of the lemma.
\endproof

Now we are ready to establish subdifferential sum rules for the restricted multifunction \eqref{restr}.

\begin{theorem}[\bf subdifferentials of restricted multifunctions]\label{indicator sum} Let $F\colon X\tto Z$ be a multifunction  between Asplund spaces, where $Z$ is ordered by the closed convex cone $\Th\subset Z$, let $\O\subset X$, and let $(\ox,\oz)\in\epi_\Th F$ with $\ox\in\O$. Assume that $F$ is epiclosed around $(\ox,\oz)$, that $\Om$ is closed around $\ox$, that the mixed subdifferential qualification condition
\begin{equation}\label{qc special sum}
{\partial}^\infty_{M,\Theta} F(x,z)\cap\big(-N(x;\Omega)\big)=\{0\}
\end{equation}
is satisfied, and that either $F$ is PSNEC at $(\ox, \oz)$ or $\Om$ is SNC at $\ox$. Then we have the sum rules
\begin{equation}\label{special sum}
\partial_\Theta F_\Om(\ox, \oz) \subset \partial_\Theta F(\ox, \oz) +N(\ox;\Om)
\end{equation}
for both subdifferentials $\partial_{N,\Theta}$ and $\partial_{M,\Theta}$, as well as the inclusion
\begin{equation}\label{special sum_sing}
\partial_\Theta^\infty F_\Om(\ox, \oz) \subset \partial_\Theta^\infty F(\ox, \oz) +N(\ox;\Om)
\end{equation}
for their singular subdifferential counterparts $\partial^\infty_{N,\Theta}$ and $\partial^\infty_{M,\Theta}$. If in addition $F$ is $N$-epiregular $($resp.\ $M$-epiregular$)$ at $(\ox, \oz)$ and $\Om$ is normally regular at $\ox$, then \eqref{special sum} holds as an equality and also the mapping $F_\O$ from \eqref{restr} has the corresponding regularity property at $(\ox, \oz)$.
\end{theorem}

{\noindent\it Proof}. Check first that $S_{\mathcal{E}}$ is inner semicontinuous at $(\ox, \oz, \oz, 0)$. Indeed, note that
$$
S_{\mathcal{E}}(x,z)=\big\{(z_1,z_2)\in Z \times Z\;\big|\;z_1\in {\mathcal{E}}_{F,\Theta}(x),\;z_2\in {\mathcal{E}}_{\Delta,\Theta}(x),\; z=z_1+z_2\big\}.
$$
Since $\oz\in{\mathcal E}_F(\ox)\subset F(\ox)+\Th$ and $0\in\Delta(\ox;\O)+\Th$ with $\ox\in\O$, it follows that $(\oz,0)\in S_{\mathcal{E}}(\ox,\oz)$, i.e., $(\ox,\oz,\oz,0)\in\gph S_{\mathcal{E}}$. To verify the inner semicontinuity of $S_{\mathcal{E}}$ at $(\ox,\oz,\oz,0)$, take any $(x_k,z_k)\to(\ox,\oz)$ from the set $\dom S_{\mathcal{E}}$ and then choose $z_{1k}:=z_k$ and $z_{2k}:=0$ for all $k\in\N$. This gives us $(z_{1k},z_{2k})\in S_{\mathcal{E}}(x_k,z_k)$ and $(z_{1k},z_{2k})\to(\oz,0)$ as $k\to\infty$. which ensures that the multifunction $S_{\mathcal{E}}$ is inner semicontinuous at $(\ox,\oz,\oz,0)$.

Next we show that the SNC property of $\Om$ at $\ox$ implies that $\epi_\Theta \Delta(\cdot;\Om)$ is PSNC at $(\ox,0)$, i.e., $\Delta(\cdot;\Om)$ is PSNEC at this point. Indeed, take arbitrary sequences $(x_k,z_k)\in \epi_\Theta \Delta(\cdot;\Om)=\Om\times \Theta$ and $(x_k^*, z_k^*)\in \Hat{N}((x_k,z_k); \Om\times \Theta)=\Hat{N}(x_k; \Om)\times \Hat{N}(z_k; \Theta)$ by the normal cone product rule. This tells us therefore that $x^*_k\in \Hat{N}(x_k; \Om)$. The imposed SNC property of $\O$ at $\ox$ ensures that the convergence $x_k^* \st{w^*}{\to} 0$ and  $\|z_k^*\|\to 0$ yields the norm convergence $\|x_k^*\|\to 0$ as $k\to\infty$. This clarifies the PSNEC property of the restricted mapping $F_\O$ in \eqref{restr} at $(\ox,\ox)$.  Similarly we can check that the latter property of $F_\O$ holds provided that $F$ is PSNEC at $(\ox,\oz)$.

Further, it follows from \eqref{3.13} that the qualification condition \eqref{qc special sum} yields the fulfillment of the mixed subdifferential qualification condition \eqref{qc_order_sum} imposed in Theorem~\ref{ordered_sum} while specified to the case of the mappings $F$ and $\Delta(\cdot;\O)$ in the sum. Applying now the sum rules \eqref{order_sum_subset} and \eqref{order_sum_sing} in Asplund spaces, we obtain the claimed results \eqref{special sum} and \eqref{special sum_sing} for all the four subdifferentials.

Next we proceed with the regularity statements of the theorem. For definiteness, consider the case where $F$ is $N$-epiregular at $(\ox,\oz)$ and $\O$ is normally regular at $\ox$, while observing that the case where $F$ is $M$-epiregular at $(\ox,\oz)$ is justified similarly. Our goal is to verify the inclusion
\begin{equation}\label{3.17}
\partial_\Theta F(\ox, \oz)+N(\ox; \Om) \subset \partial_\Theta F_\Om (\ox, \oz)
\end{equation}
under the assumptions made, which leads us to the equality in \eqref{special sum} and the $N$-epiregularity of $F_\O$ due to fulfillment of inclusion  \eqref{special sum} proved in Lemma~\ref{sub indicator}.

To furnish the proof of \eqref{3.17} and of the regularity statement for $F_\Th$ under the imposed assumptions, note first that $N(\ox;\Om)=\Hat{N}(\ox;\Om)$ by the assumed normal regularity of $\O$ at $\ox$.  We also get by the $N$-epiregularity of $F$ at $(\ox,\oz)$ that
$$
D^*{\mathcal{E}}_{F,\Theta}(\ox,\oz)(z^*)=\Hat{D}^*{\mathcal{E}}_{F, \Theta}(\ox, \oz)(z^*)\;\mbox{ whenever }\;z^* \in Z^*.
$$
To verify that $F_\Om$ is $N$-epiregular at $(\ox,\oz)$, it suffices to show that for any $x^* \in \partial_\Theta F_\Om (\ox, \oz)$ we have $x^* \in \Hat{\partial}_\Theta F_\Om (\ox, \oz)$. Indeed, it follows from the construction of $F_\Om$ and the regularity properties of $F$ and $\Om$ that $ x^*\in \Hat{\partial}_\Theta F(\ox, \oz)+\Hat {N}(\ox, \Om)$. Then we get $x ^* \in \Hat{\partial}_\Theta F_\Om(\ox, \oz)$ due to Lemma~\ref{reg sub f and indicator}. To verify \eqref{3.17}, pick any $x^*\in \partial_\Theta F(\ox, \oz)+N(\ox;\Om)$ and deduce from the equalities above that
\begin{equation*}
x^*\in \Hat{\partial}_\Theta F(\ox, \oz)+\Hat {N}(\ox;\Om)\subset \Hat{\partial}_\Theta F_\Om (\ox, \oz)= \partial_\Theta F_\Om (\ox, \oz).
\end{equation*}
Thus we arrive at \eqref{3.17} and complete the proof of the theorem.
\endproof
\vspace*{0.05in}

Observe that when both spaces $X$ and $Z$ are finite-dimensional, the PSNEC property in the results above holds automatically and the normal mixed subdifferentials are identical; the same is true for ther singular counterparts. Furthermore, the inner semicompactness property of $S_{{\mathcal{E}}}$ holds if this mapping is uniformly bounded around $(\bar x,\bar z)$.\vspace*{0.05in}

We end this section by presenting two examples illustrating both assumptions and conclusions of the obtained subdifferential sum rules. The first example confirms explicitly that the fulfillment of the sum rule \eqref{order_sum_subset} under the fulfillment of the qualification condition \eqref{qc_order_sum}.

\begin{example}[\bf illustration of the subdifferential sum rules]\label{ex-sum} {\rm Consider the two set-valued mappings $F_1,F_2\colon\mathbb R\tto\mathbb R$ given by
\begin{eqnarray*}
F_1(x):=\begin{cases}[0,\infty)&\mbox{for }\; x \le 0, \\ x, &\mbox{for }\, x > 0, \end{cases}
\end{eqnarray*}
\begin{eqnarray*}
F_2(x):=\begin{cases}-x&\mbox{for }\; x < 0, \\ [0,\infty) &\mbox{for }\, x \ge 0 \end{cases}
\end{eqnarray*}\\
with the ordering cone $\Theta:=\{x\in\mathbb R\;|\;x\ge 0\}$. Then the qualification condition \eqref{qc_order_sum} is clearly satisfied, and we have the subdifferential sum rule
\begin{eqnarray*}
\partial_{\Theta}F(0,0)= [-1,1]=\partial_{\Theta}F_1(0,0)+\partial_{\Theta}F_2(0,0),
\end{eqnarray*}
which holds in this case as an equality.}
\end{example}\vspace*{0.05in}

The next example shows that the inclusion in \eqref{order_sum_subset} may be strict.

\begin{example}[\bf strict inclusions of the subdifferential sum rules]\label{strict ex-sum} {\rm Consider the two multifunctions defined as epigraphs of the nonsmooth functions by
\begin{equation*}
F_1(x):= \left\{y\in\mathbb R\;\big|\;y \geq |x| \right\},
\end{equation*}
\begin{equation*}
F_2(x):= \left\{y \in\mathbb R\;\big|\;y \geq -|x| \right\}.
\end{equation*}
We can easily compute by the subdifferential definitions that
\begin{equation*}
\partial^\infty F_1(0,0)= \partial^\infty F_2(0,0)=\{0\},
\end{equation*}
which guarantees the fulfillment of the qualification condition \eqref{qc_order_sum}, and also
\begin{equation*}
\partial F_1 (0,0)= [-1, 1],
\end{equation*}
\begin{equation*}
\partial F_2 (0,0)= \{-1, 1\}
\end{equation*}
from which we get the equality $\partial F_1(0,0)+\partial F_2(0,0)=[-2, 2]$.
On the other hand,
\begin{equation*}
(F_1+F_2)(x)= [0, \infty)\;\mbox{ for any }\;x\in\mathbb R,
\end{equation*}
and thus $\partial (F_1+F_2)(0,0)= \{0\}$. This shows that the equality fails in \eqref{order_sum_subset}.}
\end{example}

\section{Subdifferential Chain Rules for Ordered Multifunctions}
\setcounter{equation}{0}

This section is devoted to deriving various {\em subdifferential chain rules} for general compositions of ordered set-valued mappings and their specifications.

Consider two multifunctions $G \colon X \tto Y$ and $F \colon Y \tto Z$ between normed spaces, where $Y$ and $Z$ are ordered via \eqref{order} by the closed convex cones $\Theta_1\subset Y$ and $\Theta_2\subset Z$, respectively. The {\em composition}
$(F\circ G)\colon X \tto Z$ of these mappings is defined by
\begin{equation}\label{mapping_composition}
(F\circ G)(x):=\bigcup\limits_{y\in G(x)}F(y)=\big\{z\in Z\;\big|\;\exists\,y\in G(x)\;\text{ with }\;z\in F(y)\big\}, \quad x\in X.
\end{equation}

To derive chain rules for the four subdifferentials in Definition~\ref{sub_or}(ii--v) applied to ordered compositions \eqref{mapping_composition}, we recall first the chain rules for normal and mixed coderivatives of set-valued mappings between Asplund spaces taken from  \cite[Theorem~3.13]{M1}.

\begin{lemma}[\bf chain rules for coderivatives]\label{coder_chain}
Let $G \colon X \tto Y$ and $F \colon Y \tto Z$ be multifunctions between Asplund spaces, and let $(\ox,\oz)\in\gph(F\circ G)$. Define the mapping $K\colon X\times Z\tto Y$ by
\begin{equation}\label{K}
K(x,z):= G(x)\cap F^{-1}(z)=\big\{y \in G(x)\;\big|\;z \in F(y)\big\}.
\end{equation}
Then the following assertions hold:
\begin{itemize}	
\item[{\bf (i)}] Given $\bar y \in K(\bar x, \oz) $, assume that the mapping $K$ from \eqref{K} is inner semicontinuous at $(\bar x, \oz, \bar y)$, that the sets $\gph F$ and $\gph G$ are locally closed around the point $(\bar y, \oz)$ and $(\bar x, \bar y)$, respectively, that either $F$ is PSNC at $(\bar y, \bar x)$ or $G^{-1}$ is PSNC at this point, and that the mixed coderivative qualification condition
\begin{equation}\label{qc-chain}
D_M^*F(\bar y,\oz)(0)\cap \big(-D_M^*G^{-1}(\bar x,\bar y)(0)\big)=\{0\}
\end{equation}
is satisfied. Then we have the chain rules
\begin{equation}\label{coder_chain_subset}
D^*(F\circ G)(\bar x,\oz)(z^*)\subset D_N^*G(\bar x,\bar y)\circ D^*F(\bar y,\oz)(z^*),\quad z^*\in Z^*,
\end{equation}	
valid for both normal and mixed coderivatives $D^*=D^*_N$ and $D^*=D^*_M$.

\item[{\bf (ii)}]Assume that the mapping $K$ from \eqref{K} is inner semicompact at $(\bar x, \oz)$, and that all the other conditions of {\rm(i)} hold for all $\bar y\in K(\ox,\oz)$. Then we have the inclusions
\begin{equation*}
D^*(F\circ G)(\bar x,\oz)(z^*)\subset\bigcup\limits_{\bar y\in K(\bar x,\oz)}\big[D_N^*G(\bar x,\bar y)\circ D^*F(\bar y,\oz)(z^*)\big],\quad z^*\in Z^*
\end{equation*}
valid for both normal and mixed coderivatives $D^*=D^*_N$ and $D^*=D^*_M$.
\end{itemize}	
\end{lemma}

The next lemma shows how the epigraphical multifunctions behave under compositions of ordered mappings. It is proved similarly to Lemma~\ref{epi_sumrule} in the case of summation.

\begin{lemma}[\bf epigraphical multifunctions under compositions]\label{epi- chain}
Let $G\colon X \tto Y$  and $F\colon Y \tto Z$ be set-valued mappings between normed spaces, where $Y$ and $Z$ are ordered by some closed convex cones $\Theta_1\subset Y$ and $\Theta_2\subset Z$, respectively. Then we have the equality
\begin{equation}\label{epi-chain}
\big({\mathcal{E}}_{F, \Theta_2}\circ{\mathcal{E}}_{G, \Theta_1}\big)(x)= {\mathcal{E}}_{(F\circ G, \Theta_2)}(x)\;\mbox{ for all }\;x\in X.
\end{equation}
\end{lemma}

Now we are ready to derive general chain rules for all the subdifferentials of ordered multifunctions taken from Definition~\ref{sub_or}(ii--v).

\begin{theorem}[\bf chain rules for subdifferentials of ordered multifunctions]\label{ordered_chainrule}
Let $G\colon  X\tto Y$ and $F\colon Y\tto Z$ be set-valued mappings between Asplund spaces, where $Y$ and $Z$ are ordered by the closed convex cones $\Th_1\subset Y$ and $\Th_2\subset Z$, respectively, and let $(\bar x, \oz) \in \epi_{\Theta_2}(F\circ G)$. Define the set-valued mapping $H\colon X\times Z\tto Y$ by
\begin{equation}\label{H}
H(x, z):={\mathcal{E}}_{G, \Theta_1}(x)\cap\big({\mathcal{E}}_{F, \Theta_2}\big)^{-1}(z)\;\mbox{ for all }\;(x,z)\in X\times Z.
\end{equation}
Then the following subdifferential chain rules hold:
\begin{itemize}
\item[]{\bf(i)} Fix $\bar y\in H(\bar x,\oz)$ for $H$ from \eqref{H} and suppose that this mapping is inner semicontinuous at $(\bar x,\oz,\bar y)$. Assume also that $G$ and $F$ are epiclosed around the corresponding points, that either $F$ is PSNEC at $(\bar y, \oz)$ or $G^{-1}$ is PSNEC at $(\bar y,\bar x)$, and that the qualification condition
\begin{equation}\label{qc_ordered_chain}
\partial^\infty_{M,\Theta_2}F(\bar y,\bar z)\cap\big(-\partial^\infty_{M,\Theta_2}G^{-1}(\bar y,\bar x)\big)=\{0\}
\end{equation}
is satisfied. Then we have the chain rule inclusions
\begin{equation}\label{order_chain_subset}
\partial_{\Theta_2}(F\circ G)(\ox,\oz)\subset\partial_{N,\Th_1}G(\bar x,\bar y)\circ\partial_{\Theta_2}F(\bar y,\bar z)
\end{equation}
valid for both normal and mixed subdifferentials from Definition~{\rm\ref{sub_or}(ii,iii)}. Moreover, we have the singular subdifferential chain rules
\begin{equation}\label{order_chain_sing}
\partial^\infty_{\Theta_2}(F\circ G)(\bar x,\bar z)\subset\partial^\infty_{N,\Th_1}G(\bar x,\bar y)\circ\partial^\infty_{\Theta_2}F(\bar y,\bar z)
\end{equation}
valid for both normal and mixed singular subdifferentials from Definition~{\rm\ref{sub_or}(iv,v)}.

\item []{\bf(ii)} If the inner semicontinuity assumption on the mapping $H$ in {\rm(i)}    is replaced by its semicompactness at $(\bar x,\bar z)\in\epi_{\Th_2}(F\circ G)$ and if all the conditions in {\rm(i)}  are satisfied for every $\bar y\in H(\bar x,\bar z)$, then we have the subdifferential chain rules
\begin{equation*}
\partial_{\Theta_2}(F\circ G)(\ox,\oz)\subset\bigcup\limits_{\bar y\in H(\ox,\oz)}\big[\partial_{N,\Th_1}G(\ox,\bar y)\circ\partial_{\Theta_2}F(\bar y,\oz)\big],
\end{equation*}
\begin{equation*}
\partial^\infty_{\Theta_2}(F\circ G)(\ox,\oz)\subset\bigcup\limits_{\bar y\in H(\ox,\oz)}\big[\partial^\infty_{N,\Th_1} G(\ox,\bar y)\circ\partial^\infty_{\Theta_2}F(\bar y,\oz)\big],
\end{equation*}
where the first inclusion holds for both normal and mixed subdifferentials, while the second one is fulfilled for the normal and mixed singular subdifferentials, respectively.
\end{itemize}
\end{theorem}
{\noindent\it Proof}. For definiteness, we verify only assertion (i), while observing that assertion (ii) can be proved similarly by applying Lemma~\ref{coder_chain}(ii) instead of Lemma~\ref{coder_chain}(i) below. Let us first justify inclusion \eqref{order_chain_subset} arguing in parallel for both normal and mixed subdifferentials with taking into account that products of Asplund spaces are Asplund as well; see, e.g., \cite{fabian,M1}.

Take any $x^* \in \partial_{\Theta_2}(F \circ G)(\ox, \oz)$ and by Definition~\ref{sub_or}(ii,iii) get  $x^* \in D^*({\mathcal{E}}_{F\circ G, \Theta_2})(\ox, \oz)(z^*)$ for the corresponding coderivative with $-z^* \in N(0;\Theta_2)$ and $||z^*||=1$. It is easy to deduce from the definitions that the imposed mixed subdifferential qualification condition \eqref{qc_ordered_chain} yields the mixed coderivative one \eqref{qc-chain} from Lemma~\ref{coder_chain}(i). Employing the latter result for the case of the epigraphical multifunctions ${\mathcal E}_{F,\Th_2}$ and ${\mathcal E}_{G,\Th_1}$
and using the normal and mixed subdifferential constructions again, we arrive at $x^*\in\partial_{N,\Th_1}(\bar x,\bar y)\circ\partial_{\Th_2}F(\bar y,\bar z)$, which justifies the subdifferential chain rules \eqref{order_chain_subset} for both normal and mixed subdifferentials of ordered multifunctions.

To verify the chain rules in \eqref{order_sum_sing} for both normal and mixed singular subdifferentials, we proceed similarly to the above with choosing $z^*=0$ according to Definition~\ref{sub_or}(iv,v).
\endproof
\vspace*{0.05in}

Let us discuss some issues related to Theorem~\ref{ordered_chainrule} and its proof.

\begin{remark}[\bf discussions on the subdifferential chain rules]\label{rem-chain} {\rm In the setting of Theorem~\ref{ordered_chainrule}, we have the following observations:

{\bf(i)} Although the chain rules in assertions (i) and (ii) hold for both normal and mixed subdifferentials and their singular counterparts,  only the {\em normal} subdifferential and the {\em normal singular one} are used in \eqref{order_chain_subset} and \eqref{order_chain_sing}, as well as in assertion (ii), respectively.

{\bf(ii)} On the other hand, only the {\em mixed singular} subdifferential is used in the subdifferential chain rules for all the four subdifferentials under consideration.

{\bf(iii)} Similarly to the proof of the coderivative chain rules in \cite[Theorem~3.13]{M1} (see Lemma~\ref{coder_chain} above), we may reduce the verification of Theorem~\ref{ordered_chainrule} to the setting of Theorem~\ref{indicator sum} with subdifferential sum rules for restricted multifunctions in the case where $\O:=\gph G$.}
\end{remark}

If the intermediate space $Y$ in the mapping composition \eqref{mapping_composition} is finite-dimensional, we have the subdifferential chain rules  of Theorem~\ref{ordered_chainrule} with a simplified qualification condition.

\begin{corollary}[{\bf subdifferential chain rules with finite-dimensional intermediate spaces}]\label{cor-inter}\,
Assume that $\dim Y<\infty$ in the setting of Theorem~{\rm\ref{ordered_chainrule}}. Then all the statements of the theorem hold with the replacement of the qualification condition \eqref{qc_ordered_chain} by the following one:
\begin{equation}\label{qc chain finite}
\partial^\infty_{M,\Theta_2}F(\bar y,\oz)\cap\ker\partial^\infty_{\Theta_1}G(\bar x,\bar y)=\{0\},
\end{equation}
where the latter subdifferential stands for the common subdifferential $\partial^\infty_{\Theta_2}=\partial^\infty_{N,\Theta_2}=\partial^\infty_{M.\Theta_2}$.
\end{corollary}
{\noindent\it Proof}. We have the equality $\partial^\infty_{N,\Theta_1}G(\bar x,\bar y)=\partial^\infty_{M.\Theta_1}G(\bar x,\bar y)$ in \eqref{qc chain finite} due to the singular subdifferential definitions and the fact that there is no difference between normal and mixed coderivatives of mappings with finite-dimensional range spaces. Moreover, we can deduce from the definitions the equivalence
\begin{equation*}
z \in D^*H(\ox, \bar y)(u) \Longleftrightarrow -u \in D^*H^{-1}(\bar y, \ox)(-z)
\end{equation*}
for any multifunction $H\colon X\tto\mathbb{R}^n$. This implies that the qualification  condition \eqref{qc_ordered_chain} reduces  to \eqref{qc chain finite} provided that the space $Y$ is finite-dimensional.
\endproof
\vspace*{0.05in}

The next corollary is a specification of Theorem~\ref{ordered_chainrule} in the case of single-valued inner mappings in composition \eqref{mapping_composition} and reveals some conditions ensuring equalities in the chain rules.

\begin{corollary}[\bf subdifferential chain rules as equalities]\label{chain single}
Let $G=g$ be single-valued and Lipschitz continuous around $\ox$, which implies that $H$ in \eqref{H} is inner semicompact at $(\ox, \bar z)$. Assume in addition that in Theorem~{\rm\ref{ordered_chainrule}(ii)}, the multifunction $F$ is $N$-epiregular $($resp.\ $M$-epiregular$)$ at $(\bar y,\bar z)$ with $\bar y:= g(\ox)$, and that either $g$ is $N$-regular at $\ox$ with $\dim Y < \infty$, or $g$ is strictly differentiable at $\ox$. Then the composition $F\circ g$ is $N$-epiregular at $(\ox, \bar z)$ and we have the following chain rule equality for both normal and mixed subdifferentials of the composition:
\begin{equation*}
\partial_{\Theta_2}(F\circ g)(\ox,\bar z) =\partial_{N,\Th_1}g(\ox)\circ\partial_{\Theta_2}F(\bar y, \bar z).
\end{equation*}
\end{corollary}
{\noindent\it Proof}. Follows from \cite[Theorem~3.13]{M1} and the subdifferential definitions.
\endproof

Finally in this section, we present effective sufficient conditions that ensure the fulfillment of the mixed subdifferential qualification condition \eqref{qc_ordered_chain} together with the PSNEC requirements, and hence of the subdifferential chain rules in Theorem~\ref{ordered_chainrule} under the {\em well-posedness} properties of the ordered multifunctions involved into compositions. Here we benefits from the {\em coderivative criterion} for the Lipschitz-like property presented in Theorem~\ref{lipschitz like} and its equivalent form characterizing metric regularity as discussed above.

\begin{corollary}[\bf subdifferential chain rules for ELL and metrically regular mappings]\label{chain-mr} $\,$ Suppose that all the assumptions in Theorem~{\rm\ref{ordered_chainrule}(i)} but the qualification condition \eqref{qc-chain} and the PSNEC properties imposed therein, are satisfied. Assume also that either $F$ is ELL around $(\bar y, \bar z)$, or $G$ is metrically regular around $(\ox, \bar y)$. Then the subdifferential chain rules of Theorem~{\rm\ref{ordered_chainrule}(i)} hold for both mixed and normal subdifferentials in \eqref{order_chain_subset} and their singular counterparts in \eqref{order_chain_sing}.
\end{corollary}
{\noindent\it Proof}. As follows from the applications of Theorem~\ref{lipschitz like} and its equivalent version for metric regularity (see \cite[Theorem~1.49]{M1}) to the case of epigraphical multifunctions, the well-posedness conditions imposed in the corollary yield the simultaneous fulfillment of the qualification condition \eqref{qc_ordered_chain} and the PSNEC properties assumed in Theorem~\ref{ordered_chainrule}(i). This completes the proof.
\endproof

\section{Concluding Remarks}
\setcounter{equation}{0}

In this paper, we defined several subdifferential constructions for set-valued mappings between normed spaces and then established fundamental sum and chain rules for such subdifferentials in the framework of Asplund spaces, which constitute a broad class of Banach spaces including, in particular, all reflexive ones. The obtained calculus rules are the first results of subdifferential calculus for ordered multifunctions between infinite-dimensional spaces.

Our future goals are to apply the developed subdifferential calculus rules to establishing existence theorems and deriving necessary optimality conditions in problems of vector and set-valued optimization (including equilibrium models) with structured  constraints of various types: functional, operator, equilibrium ones, etc. Some results in this direction were obtained in finite-dimensional spaces in our previous publication \cite{MO}, while the infinite-dimensional settings are more challenging and cover a large variety of mathematical models and practical applications.\\[1ex]
{\bf Acknowledgements}. The authors are grateful to two anonymous referees for their careful reading the paper with providing important remarks and references, which allowed us to improve the original presentation.

\end{document}